\documentclass[11pt]{article}
\usepackage{amssymb}
\usepackage{graphicx}
\usepackage{amsmath}
\usepackage{xcolor}
\usepackage{colortbl,booktabs}
\usepackage{multirow}
\usepackage{subfigure}
\usepackage{multicol}
\usepackage{cite}
\usepackage{epstopdf}

\def\R{{\mathbb R}} 
 \def\N{{\mathbb N}}

\newenvironment{prof}[1][Proof]{\noindent\textit{#1}\quad }
{\hfill $\Box$\vspace{0.7mm}}
 \def\qq{\qquad}

\def\bb{\begin}
\def\bc{\begin{center}}
\def\ec{\end{center}}
\def\ba{\begin{array}}
\def\ea{\end{array}}

\def\be{\begin{equation}}     \def\ee{\end{equation}}
\def\bea{\begin{eqnarray}}    \def\eea{\end{eqnarray}}
\def\beaa{\begin{eqnarray*}}  \def\eeaa{\end{eqnarray*}}

\def\hh{\!\!\!\!}             
\def\EQ{\hh & = & \hh}

              \def\vp{\varphi}

\def\nn{\nonumber}            
\def\oo{\infty}               
               \def\pa{\partial}
                 \def\qq{\qquad}
\def\f{\frac}                 
\def\z{\left}                 \def\y{\right}

\def\rd{\,{\rm d}}
\def\dt{\,{\rm d}t}

\def\du{\,{\rm d}u}

\def\ifl{\iffalse}
\def\qqf{\qquad \forall}
\def\andq{\quad \mbox{ and } \quad}

\def\lb{\label}
\def\x#1{(\ref{#1})}

\def\pas#1#2{\left.#1\right|_{#2}}

\def\pp{{m,p}}
\def\vpm{\eta_\pp}
\def\T{{2m\pi}}
\def\imt{\int_0^\T}

\def\ol{\overline}

\def\ol{\overline}
\def\cc{C}
\def\ss{S}

\begin{document}

\begin{center}
{\LARGE On the Stability of Symmetric Periodic Orbits of a Comb-Drive Finger Actuator Model }\\
\vskip 0.3cm

Xuhua Cheng\footnote{Corresponding author. This author is supported by the National Natural Science Foundation of China (Grant No. 11601257) and Natural Science Foundation of Hebei Province (Grant No. A2019202342).} \\
School of Science, Hebei University of Technology, \\ Tianjin 300130, China\\
E-mail: {\tt 2017081@hebut.edu.cn}
\vskip 0.2cm

Baoting Liu \\
School of Science, Hebei University of Technology, \\ Tianjin 300130, China\\
E-mail: {\tt liubaoting2022@163.com}

\end{center}
\vskip 0.2cm

\bb{abstract}

In this paper, we study the stability of symmetric periodic solutions of the comb-drive finger actuator model.
First, on the basis of the relationship between the potential and the period as a function of the energy, we derive the properties of the period of the solution of the corresponding autonomous system (the parameter $\delta$ of input voltage $V_\delta(t)$  is equal to zero) in the prescribed energy range. Then, using these properties and the stability criteria of symmetric periodic solutions of the time-periodic Newtonian equation, we analytically prove the linear stability/instability of the symmetric $(m,p)$-periodic solutions which emanated from nonconstant periodic solutions of the corresponding autonomous system when the parameter $\delta$ is small.


\end{abstract}

{\bf Mathematics Subject Classification (2010):} 70H14; 34D20; 34C25;

{\bf Keywords:} comb-drive finger actuator model; stability; symmetric periodic solutions; hyperbolic

\renewcommand{\raggedright}{\leftskip=0pt \rightskip=0pt plus 0cm}

\section{Introduction}


Whether in mathematics or anything else, periodic solutions of conservative system have always been an important research field, which has attracted attention of many researchers. We just refer the reader to \cite{HLS14,LO08,LS80,M1981,O96,J1979,RP1982,WR1970,AF1995,ML2014,AM1994,SZ1990,SZ1997,FJ2003} and the references. However, in terms of the analysis of stability of the periodic solutions, it is less explored and most known stability results are based on numerical calculations or concerned with the stability of the equilibriums. Moreover, for conservative systems, the classics Lyapunov's methods do not more for studying the Lyapunov stability/instability. So, it is involved of deep theory like the Birkhoff normal forms and Moser's twist theorem, even for systems of lower degrees of freedom \cite{M1, SM}. In fact, even for the equilibria of non-autonomous systems, it is also rather difficult, much more nonconstant periodic solutions.

As a beginning step towards the stability of nonconstant periodic solutions, Zhang and his co-authors established stability criteria of nonconstant symmetric periodic solutions of a second-order nonlinear scalar Newtonian equation based on the theory for Hill equation in \cite{CC,MW,ZCC18}, where the symmetries mean the oddness or the evenness of the periodic solutions in time. Applying the stability criteria, some analytical results were derived on the linear stability/instability of odd and even $(m,p)$-symmetric periodic solutions of the elliptic Sitnikov problem \cite{CC,CL,Ch}.

Comb-drive finger actuator, as a special type of micro-electromechanical-system (MEMS) \cite{TV:2006,SJ:2007,F:1991,D:2005,MT:2008,MIY2011,LZ2022,Y:2012,CL2010}, is widely used in sensing and actuation, such as resonant sensors \cite{TV:2006}, accelerometers \cite{F:1991} and optical communication devices\cite{MT:2008}. An excellent bibliographic source of applications of this devices can be founded in \cite{MIY2011} and references therein.
In recent years, there are some results on periodic solutions of this model. For example,
 Guti\'{e}rrez et al. firstly studied the existence and stability of periodic oscillations in canonical MEMS like \cite{HC:1967} by using topological and variational methods~\cite{AP:2013}. Then, Alexander et al. studied the existence of even and $nT$-periodic solutions ($n$ positive integer) for the MEMS in a quantifiable $\delta$-interval by using global continuation method of the zeros of a function depending on one parameter provided by the Leray and Schauder Theorem
 \cite{ADA:2017}. Besides, N\'{u}\~{n}ez et al. presented the existence of odd and $nT$-periodic solutions ($n$ is a positive integer) for the transverse comb-drive finger device using shooting method and the Sturm comparison theory~\cite{DOL:2021}.


However, the literature \cite {ADA:2017} indicated that stability of the nonconstant even symmetric periodic solutions of comb-drive finger actuator model is not easy issue.
Motivated by the works of the literature~\cite{CC,CL,ZCC18,Ch}, we in this paper try to study the stability/instability of these nonconstant symmetric periodic solutions of the literature~\cite{ADA:2017,DOL:2021} by stability criteria established in \cite{ZCC18}. Specifically, based on the relationship between the potential and the period as a function of the energy, we firstly deduce the properties of the period of the solution of the corresponding autonomous equation (the parameter $\delta$ of input voltage $V_\delta(t)$  is equal to zero) in the prescribed energy range. Then, according to these properties and stability criteria of \cite {ZCC18}, we analytically prove the linear stability/instability of the symmetric $(m,p)$-periodic solutions which emanated from nonconstant periodic solutions of the corresponding autonomous equation when the parameter $\delta$ is small. To our knowledge, it is the first time to analytically study the stability of nonconstant periodic solutions for comb-drive finger actuator model. Meanwhile, it further reveals the mathematical mechanism of the complex dynamics in this model.

The rest of the paper is organized as follows. In Section 2, we will recall stability criteria of symmetric periodic solutions of a second-order nonlinear scalar Newtonian equation in \cite {ZCC18}. In Section 3, comb-drive finger actuator model will be introduced and some properties of the period of the solution of the corresponding autonomous equation will be derived. In Section 4, we will arrive at some results on the linear stability/instability of nonconstant symmetric periodic solutions. Finally, we will give a conclusion.


\section{Preliminaries} \setcounter{equation}{0} \lb{pre}

In this section, we will briefly recall the stability criteria of symmetric periodic solutions in a second-order nonlinear scalar Newtonian equation~\cite{ZCC18,CC}.

Consider the following second-order nonlinear scalar Newtonian equation
    \be \lb{xe}
    \ddot x+F(x,t,e)=0,
    \ee
where $0\leq e\leq1$ and $F(x,t,e)$ is a smooth function of $(x,t,e)\in \R^3$ and satisfies the symmetries of \x{Sy1} and the minimal period in $t$ is $T$.
    \be \lb{Sy1}
    \z\{\ba{l} F(-x,t,e) \equiv -F(x,t,e), \\
    F(x,-t,e) \equiv F(x,t,e),\\
    F(x,t+2\pi,e) \equiv F(x,t,e),\\
    F(x,t,0)\equiv f(x) \mbox{ and } f(x) \mbox{ is }  \mbox{ odd },  \\
    x f(x)> 0\mbox{ for }x\ne 0.
    \ea\y.
    \ee
when $e=0$, it corresponds to the following autonomous system:
    \be \lb{x}
    \ddot x+f(x)=0,
    \ee

Because of the symmetries of $f(x)$, there are the following two classes of periodic solutions of Eq.\x{x}: odd periodic solutions and even periodic solutions.

\bb{Definition}\lb{D1}
Let $x(t)=\ss(t)=\ss(t,\eta)$ be the solution of Eq. \x{x} satisfying the initial value conditions
    \be \lb{ini}
    \z(x(0), \dot x(0)\y)=(0,\eta),
    \ee
where $\eta\in \z(0,\eta_{\max}\y)$ and $\eta_{\max}:=\sqrt{2 E_{\max}}$.
If $\ss(t)$ is a periodic solution of Eq.\x{x} with the minimal period
    \be \lb{Tv}
    T=T(h),\qq \mbox{where } h=\eta^2/2,
    \ee
and satisfy with the following symmetries
    \be \lb{Os1}
    \ss(-t) \equiv -\ss(t) \andq \ss(t+T/2)\equiv - \ss(t),
    \ee
then $\ss(t)$ is called  an odd periodic solution of Eq. \x{x}.
\end{Definition}


\bb{Definition}\lb{D2}
Let $x(t)=\cc(t)=\cc(t,\xi)$ be the solution of Eq. \x{x} satisfying the initial value conditions
    \be \lb{ini2}
    \z(x(0), \dot x(0)\y)=(\xi,0),
    \ee
where $\xi\in \z(0,+\oo\y)$. If $\cc(t)$ is a periodic solution of Eq. \x{x} with the minimal period
    \be\lb{Tv1}
    T=T(h),\qq \mbox{where } h=E(\xi),
    \ee
and satisfy with the following symmetries
    \be \lb{Os2}
    \cc(-t)\equiv \cc(t) \andq \cc(t+T/2)\equiv - \cc(t),
    \ee
then $\cc(t)$ is called  an even periodic solution of Eq. \x{x}.
\end{Definition}

 Notice that $\ss(t)>0$ is strictly increasing and $\cc(t)>0$ is strictly decreasing on $(0,T/4)$, respectively. Moreover,
the solutions $\ss(t)$ and $\cc(t)$ are $T/2$-anti-periodic from \x{Os1} and \x{Os2}, i.e,
    \be \lb{Os12}
    \ss(T/2-t) \equiv \ss(t) \andq \cc(T/2-t) \equiv -\cc(t).
    \ee
In addition, if $\eta$ and $\xi$ satisfy
    \be \lb{sc1}
     \eta^2/2=E(\xi)=:h,
    \ee
then $\ss(t)$ and $\cc(t)$ satisfy
    \be \lb{sc2}
    \ss(t+T/4) \equiv \cc(t) \andq \cc(t+T/4) \equiv -\ss(t).
    \ee

\bb{Definition}\lb{D3}
 Let $x(t)$ be an $(\pp)$-periodic solution of Eq. \x{x}, if $x(t)$ is a $\T$-periodic solution of Eq. \x{x} and has precisely $2p$ zeros in intervals $[t_0,t_0+\T)$, where $m,~p\in \N$.
\end{Definition}

Let $\phi_\pp(t):= \ss(t,\vpm)$ be an $(\pp)$-odd periodic solution of Eq. \x{x} with the minimal period $T(h_{\pp})= \T/p$, where $\vpm: = \sqrt{2 h_\pp}$.
Then when $T'(h_\pp)\ne 0$, there exists an odd $(\pp)$-periodic solution $\phi_\pp(t,e)$ of Eq. \x{xe} emanating from $\phi_\pp(t)$ (see \cite{CC}, Theorem 3.1). Similarly,
Assuming that $\xi_\pp>0$ and $E(\xi_\pp)=h_\pp$,
let $\vp_\pp(t):= C(t,\xi_\pp)$
be an even $(\pp)$-periodic solution of \x{x} of the minimal period $T(h_{\pp})=\T/p$.
So, there also exists an even $(\pp)$-periodic solution $\vp_\pp(t,e)$ of Eq. \x{xe} emanating from $\vp_\pp(t)$ when $T'(h_\pp)\ne 0$ (see \cite{CC}). Note that $\phi_\pp(t)$ and $\vp_\pp(t)$ are $\T/p$-periodic. Usually speaking, if $e>0$, the minimal period of $\phi_\pp(t,e)$ and $\vp_\pp(t,e)$ are $\T$, not $\T/p$.

We know that $\phi_\pp(t,e)$ is the odd $(\pp)$-periodic solution of Eq. \x{xe}. Then, for arbitrary $e\in[0,e_\pp)$, the linearization equation of Eq. \x{xe} along $x=\phi_\pp(t,e)$ is the following Hill equation
     \be \lb{He}
    \ddot y + q(t,e)y=0, \qq q(t,e):=\pas{\f{\pa F}{\pa x}}{(\phi_\pp(t,e),t,e)}.
    \ee
Assuming that the period of $\phi_\pp(t,e)$ is $T=\T$, the corresponding trace of the $2m\pi$-periodic Poincar$\acute{\hbox{e}}$ matrix of Eq. \x{He} is
    \be \lb{Tre}
\ol{\tau}_\pp(e):={\ol{\psi}_1(\T,e)}+{\dot{\ol{\psi}}_2(\T,e)},
    \ee
where $\ol{\psi}_i(t,e)$ $(i=1,2)$ are fundamental solutions of Eq. \x{He}. Based on theory for the Hill equation \cite{MW}, the stability criteria of the odd $(m,p)$-periodic solution $\phi_\pp(t,e)$ are given as follows.

 \bb{Lemma}~(\cite{CC}, Theorem 3.4 and Corollary 3.5) \lb{The1}
Let $\phi_\pp(t)$ be the odd $(\pp)$-periodic solution of Eq. \x{x} verifying condition $T'(h_\pp)\ne 0$. Denote
    \be \lb{F23}
    \ol{F}_{23}(t):=\pas{\f{\pa^2 F}{\pa t\pa e }}{(\phi_\pp(t),t,0)}.
    \ee
Then the derivative of the trace of (\ref{Tre}) at $e=0$ is
    \be \lb{dtau0}
    \ol{\tau}_\pp'(0):= \pas{\f{\rd\ol{\tau}_\pp(e)}{\rd e}}{e=0} =-p T'(h_\pp) \imt \ol{F}_{23}(t)\dot\phi_\pp(t)\dt,
    \ee
where $h_\pp=\eta_\pp^2/2$.
Moreover, if~~$\ol{\tau}_\pp(0)=2$, then

{\rm (i)} when $\tau'_\pp(0)<0$, $\phi_\pp(t,e)$ is elliptic and is linearly stable for $0<e\ll 1$.

{\rm (ii)} when $\tau'_\pp(0)>0$, $\phi_\pp(t,e)$ is hyperbolic and is Lyapunov unstable for $0<e\ll 1$.

\end{Lemma}

\bb{Remark} \lb{Ev3.8}

{\rm In Lemma \ref{The1}, the same is true replacing ``odd'' by ``even'', ``$\phi_\pp(t)$" by ``$\varphi_\pp(t)$" and ``$\phi_\pp(t,e)$" by ``$\varphi_\pp(t,e)$".
Here $h_\pp=E(\xi_\pp)$, $e\in[0,\tilde{e}_\pp)$, ``$\varphi_\pp(t)$" and ``$\varphi_\pp(t,e)$" are the even $(\pp)$-periodic solutions of Eq. \x{x} and Eq. \x{xe}, respectively. Please see Theorem 3.1 in \cite{ZCC18}}.

\end{Remark}

\section{ The comb-drive finger actuator}\setcounter{equation}{0}\lb{CD}

In this section, we will introduce the comb-drive finger actuator model \cite{DOL:2021}, which satisfies the symmetric properties of Eq. \x{xe} in section \ref{pre}.
\label{sec:2}

The comb-drive finger actuator is a special model of the micro-electromechanical systems (MEMS). We describe it as follows: a moveable electrode (finger) is sandwiched between two stationary electrodes and moves in the transverse direction to the longitudinal axis of the stationary electrodes. Besides, the moveable finger with mass $m$ is attached to a linear spring with stiffness coefficient $k>0$ and is at the center of the two stationary electrodes at a distance $d$ (see Fig.1).
\begin{figure}
\centering
\includegraphics[height=5cm, width=7cm]{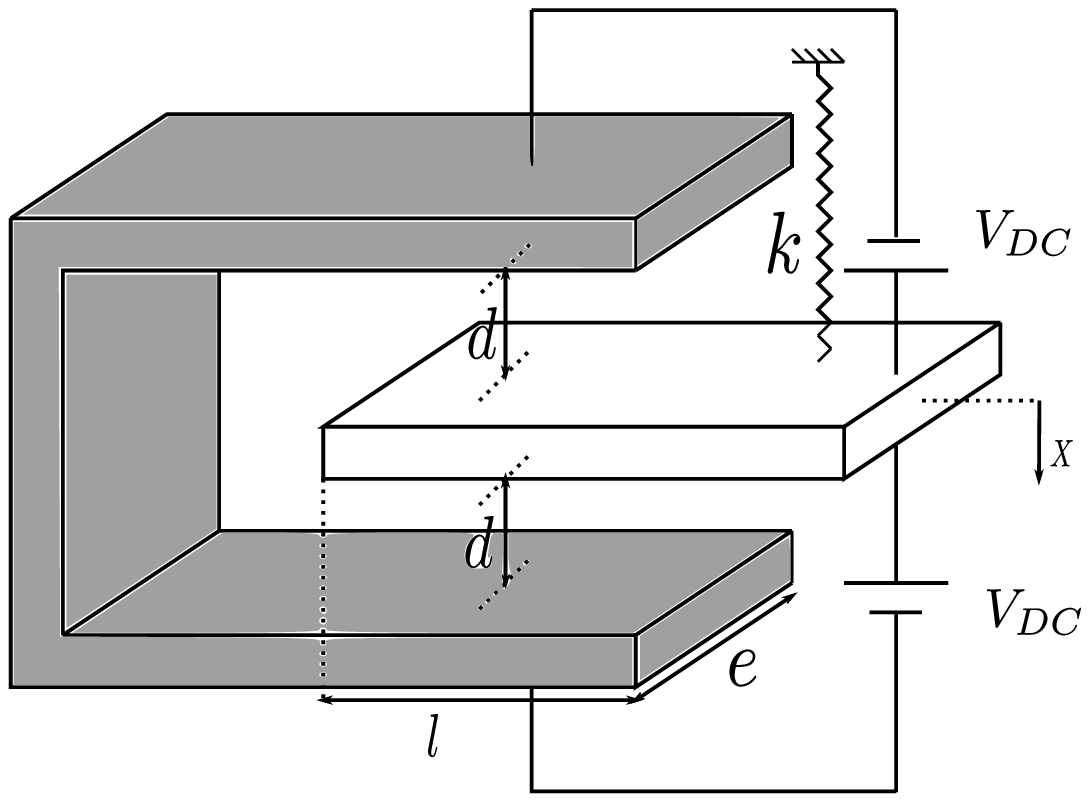}
\centerline{Fig.1  ~Idealized model of a Comb-drive finger actuator}
\end{figure}
The equation of motion of the moveable finger in a comb-drive actuator is given by
\begin{equation}\label{E21}
\ddot x+\omega^2x=\frac{4dhxV^2_\delta(t)}{(d^2-x^2)^2},~~~~~|x|<d,
\end{equation}
where $x$ denotes the vertical displacement of the movable finger from its rest position, $\omega^2=k/m$, $h=el\varepsilon/2m >0$, $e$ is the width of the plates, $l$ is the length of the fingers in the interacting zone, $\varepsilon$ is the dielectric constant of vacuum, and $V_\delta(t)$ is the input voltage. If we choose the appropriate unit of the distance and time, the Eq. \x{E21} is rewritten as
\begin{equation}\label{E211}
\ddot x+x\z(1-\frac{4\beta V^2_\delta(t)}{(1-x^2)^2}\y)=0,~~~~~|x|<1,
\end{equation}
where $\beta>0$ is a physical constant given by $\beta=\frac{el\varepsilon}{2kd^{3}}$.
Moreover, we will consider a $T_{v}$-periodic $DC$-$AC$ voltage $V_\delta(t)$ of the form
\begin{equation}\label{E22}
V_\delta(t)=V_0+\delta P(t),
\end{equation}
where $V_0 > 0$, $\delta > 0$, and $P(t)\in C(\mathbb{R}/T\mathbb{Z})$ is an even function such that $\int _{0}^{T_{v}}P(t)dt=0$. It is typically that $P(t)=cos(\omega_0 t)$ with $\omega_0=\frac{2\pi}{T_{v}}$. Letting $V_\delta(t)> 0$, we have $\delta\in[0, ~\Delta_0]$, where $\Delta_0=[0, ~\frac{-V_0}{P_m}]$ and $P_m=\min\limits_{t\in\mathbb{R}}P(t)< 0$.

Furthermore, Eq. \x{E211} is reformulated as follows:
 \be \label{E23}
\ddot x +F(x,t,\delta) =0,\qq F(x,t,\delta):= x\z(1-\frac{4\beta V^2_\delta(t)}{(1-x^2)^2}\y),
    \ee
where $F(x,t,\delta)$ is a smooth function of $(x,t,\delta)\in (-1,1)\times\mathbb{R}\times\mathbb{R}$ and
satisfies with the following symmetries
    \be \lb{Sy2}
    \z\{\ba{l} F(-x,t,\delta) \equiv -F(x,t,\delta), \\
    F(x,-t,\delta) \equiv F(x,t,\delta),\\
    F(x,t+T_{v},\delta) \equiv F(x,t,\delta),\\
    F(x,t,0)\equiv f(x) \mbox{ and } f(x) \mbox{ is }  \mbox{ odd },  \\
    x f(x)> 0\mbox{ for }x\ne 0.
    \ea\y.
    \ee
We note that $F(x,t,\delta)$ satisfies the symmetries of \x{Sy1} and the minimal period in $t$ is $T_{v}$.

In particular, when $\delta=0$, Eq. (\ref{E21}) is the following autonomous equation
\begin{equation}\label{E31}
\ddot x+x\z(1-\frac{4\beta V^2_0}{(1-x^2)^2}\y)=0.
\end{equation}
Let $V^\ast=\frac{1}{2\sqrt{\beta}}$. If $0< V_0< V^\ast$, then Eq. (\ref{E31}) has three equilibria in $(-1, 1)$ given by $(0,0)$, $(x_\ast, 0)$, $(-x_\ast, 0)$, where
$x_\ast=\sqrt{1-2V_0\sqrt{\beta}}$.
Otherwise, the moveable finger will get stuck to one of the other fixed electrodes when $V_0 > V^\ast$. This effect is known in the comb-drive literature as pull-in of double-sided capacitors or side instability.

Eq. \x{E31} is reformulated as follows:
 \be \lb{E32}
\ddot x +f(x) =0,\qq f(x):= x\z(1-\frac{4\beta V^2_0}{(1-x^2)^2}\y),
    \ee
Let the energy be $E(x)$ and $E(x)$ is formulated as follows:
    \[
    E(x)=\int_0^x f(u) \du = \frac{1}{2}x^2-\frac{2\beta V^2_0}{1-x^2}+2\beta V^2_0,
    \]
Clearly, $E(x)$ is an even function with $E(0)=0$ and $E(x)>0$ for $x\ne 0$. Then, the energy levels of solutions $x(t)$ of Eq. \x{E32} are
    \be \lb{H}
   \Gamma_{\hbar}: H(x, \dot x)=\frac{\dot x^2}{2}+\frac{1 }{2}x^2-\frac{2\beta V^2_0}{1-x^2}+2\beta V^2_0= \hbar,
    \ee
where
     \(
     \hbar\in(-\infty, \hbar_\ast]
     \)
and $\hbar_\ast:=H(x_\ast, 0)$. For $\hbar=0$,  $\Gamma$ is just the equilibrium $O(0,0)$ of Eq. \x{E32}. For $\hbar\in(0, \hbar_\ast]$,  $\Gamma_{\hbar}$ corresponds to a periodic orbit of \x{E32} with minimal period denoted by $T(\hbar)$. Next, we will give three properties of $T(\hbar)$ by the relationship between the potential and the period as a function of the energy in the following theorem.

\bb{Theorem}\label{Theo3.1}  The period $T(\hbar)$ has the following three properties:

   {\rm (i)} $\lim_{\hbar\to 0+} T(\hbar) = 2\pi/\sqrt{1-4\beta V^2_0};$

   {\rm (ii)} $\lim_{\hbar\to\hbar_\ast-} T(\hbar) =+\oo;$

   {\rm (iii)} $T'(\hbar)=\f{\rd T(\hbar)}{\rd \hbar} >0, \qqf \hbar\in (0, \hbar_\ast).$

    \end{Theorem}

In order to prove Theorem \ref{Theo3.1}, we need to introduce the following lemma.

\bb{Lemma}~(\cite{L91}, Formula(1.5) and Appendix A3) \lb{L32}
Consider the Lagrange equation
\be\lb{L}
 \ddot x + V'(x)=0,
 \ee
where $V(x)$ is a smooth potential such that
\be\lb{V}
 min_{x} V(x)=V(x_{0})=0 \qq and \qq V'(x)\neq0~for~x\neq x_{0}.
 \ee
Then, for any $h>0$, the level curve of energy $h$
$$\frac{1}{2}\dot x^{2}+V(x)=h$$
is a periodic orbit of the minimal period $T(h)$ and
\be\lb{T}
   T(h)=\sqrt{2}\int^{x_{+}(h)}_{x_{-}(h)}\frac{dx}{\sqrt{h-V(x)}},
\ee
where $x_{-}(h)$, $x_{+}(h)$ are solutions of $V(x)=h$ and satisfy $x_{-}(h)<x_{0}<x_{+}(h)$.
Moreover, the derivative of $T(h)$ in $h$ is given by
\be\lb{dT}
   \frac{dT(h)}{dh}=\frac{\sqrt{2}}{2h}\int^{x_{+}(h)}_{x_{-}(h)}\z(1-\frac{2V(x)V''(x)}{(V'(x))^{2}}\y)
   \frac{dx}{\sqrt{h-V(x)}}.
\ee
\end{Lemma}

Successively, we will prove Theorem \ref{Theo3.1} by Lemma \ref{L32}.

 \begin{prof}[~$\mathbf{Proof~of~Theorem~\ref{Theo3.1}}$]
 Assume that $x(t)$ is the $T(\hbar)$-periodic solution of Eq. \x{E32} satisfying the initial value $x(0)=0$, $\dot x(0)=\eta >0$, and $x_{-}(\hbar)$, $x_{+}(\hbar)$ are two solutions of $V(x)=E(x)=\hbar$. Then, based on  the formula \x{T} and symmetry of Eq. \x{E32}, we have
\bea \lb{T1}
  T(\hbar)\EQ \sqrt{2}\int^{x_{+}(\hbar)}_{x_{-}(\hbar)}\z(\hbar-V(x)\y)^{-\frac{1}{2}}dx\nn\\
  \EQ 2\sqrt{2}\int^{x_{+}(\hbar)}_{0}\z(\hbar-\frac{1 }{2}x^2+\frac{2\beta V^2_0}{1-x^2}-2\beta V^2_0\y)^{-\frac{1}{2}}dx\nn\\
  \EQ 2\sqrt{2}\int^{x_{+}(\hbar)}_{0}\z(\frac{1 }{2}(x^2_{+}(\hbar)-x^2)-\frac{2\beta V^2_0}{1-x^2_{+}(\hbar)}+\frac{2\beta V^2_0}{1-x^2}\y)^{-\frac{1}{2}}dx \nn\\
  \EQ  2\sqrt{2}\int^{x_{+}(\hbar)}_{0}(x^2_{+}(\hbar)-x^2)^{-\frac{1}{2}}f(x, x_{+}(\hbar))dx \nn\\
  \EQ
  2\sqrt{2}\int^{1}_{0}(1-\rho^2)^{-\frac{1}{2}}f(\rho x_{+}(\hbar), x_{+}(\hbar))d\rho,
\eea
where $\rho=\frac{x}{x_{+}(\hbar)}$ and
\bea \lb{f1}
f(x,x_{+}(\hbar))=\z(\frac{1}{2}-\frac{2\beta V^2_0}{(1-x^2)(1-x^2_+(\hbar))}\y)^{-\frac{1}{2}}.
\eea
Notice that $$\lim_{x_{+}(\hbar)\rightarrow0}f(x,x_{+}(\hbar))=\lim_{x_{+}(\hbar)\rightarrow0}f(\rho x_{+}(\hbar), x_{+}(\hbar))=\frac{\sqrt{2}}{\sqrt{1-4\beta V^2_0}}$$
and $$\lim_{x_{+}(\hbar)\rightarrow x_{\ast}}f(x,x_{+}(\hbar))=\lim_{x_{+}(\hbar)\rightarrow x_{\ast}}f(\rho x_{+}(\hbar), x_{+}(\hbar))=+\infty.$$
So, based on the dominated convergence theorem and monotone convergence theorem, we obtain
$$\lim_{\hbar\rightarrow{0}^+}T(\hbar)=\lim_{x_{+}(\hbar)\rightarrow0}T(\hbar)=\frac{2\pi}{\sqrt{1-4\beta V^2_0}}$$
and $$\lim_{h\rightarrow h_{\ast}^-}T(\hbar)=\lim_{x_{+}(\hbar)\rightarrow x_{\ast}}T(\hbar)=+\infty.$$

In the following, we continue to prove the above item (iii) of Theorem \ref{Theo3.1}. From \x{dT}, we have
\begin{align} \lb{dT1}
\frac{dT(\hbar)}{d\hbar}= &\frac{\sqrt{2}}{2\hbar}
\int^{x_{+}(\hbar)}_{x_{-}(\hbar)}
\z(\frac{4\beta V^2_0x^4[(1-x^2)(3+x^2)-12\beta V^2_0]}{x^2[(1-x^2)^2-4\beta V_0^2]^2}\y)\nonumber\\&
\z(\frac{1 }{2}(x^2_{+}(\hbar)-x^2)-\frac{2\beta V^2_0}{1-x^2_{+}(\hbar)}+\frac{2\beta V^2_0}{1-x^2}\y)^{-\frac{1}{2}}dx.
\end{align}
Let $$v(x)=(1-x^2)(3+x^2)-12\beta V^2_0.$$
Since $1-x^2>4\beta V^2_0$, we have $v(x)>0$, implying that the integrand of \x{dT1} is positive. Hence, $T'(\hbar)=\frac{dT(\hbar)}{d\hbar}>0$ for $\forall \hbar\in(0, \hbar_{*})$.
\end{prof}

\bb{Remark}\lb{R51}
{\rm Notice that the origin is surrounded by a family of periodic orbits, whose minimal period is $T(\hbar_\pp)=mT_{v}/p\in (2\pi/\sqrt{1-4\beta V^2_0},+\oo)$, i.e.,
the integers $m, \ p$ satisfy
    \be \lb{mp}
    1\le p \le \nu_m := [mT_{v}\sqrt{1-4\beta V^2_0}/2\pi],\qq m\in \N.
    \ee
In addition, the condition $T'(\hbar_\pp)\ne 0$ is ensured by (iii) in Theorem \ref{Theo3.1}.
}
\end{Remark}

\section{Main results}\setcounter{equation}{0}\lb{res}
\label{sec:4}
In the section, we will use stability criteria described in Section \ref{CD} to prove the stability/ instability of the symmetric periodic solutions of the comb-drive finger actuator model. A rigorous mathematical study on the existence and continuation of nonconstant symmetric periodic solutions of \x{E23} can be found in the literatures \cite{ADA:2017,DOL:2021}. In the following, we use the families $\phi_\pp(t,\delta)$ and $\vp_\pp(t,\delta)$ $(0 < \delta\ll \Delta_0)$ to denote the odd and even $(\pp)$-periodic solutions in $t$ of Eq. \x{E23}, respectively.

\subsection{stability of odd periodic solutions}

In this subsection, we will consider stability/ instability of the family $\phi_\pp(t,\delta)$ $(0 < \delta\ll \Delta_0)$ of odd $(\pp)$-periodic solutions of Eq. \x{E23} for $m, \ p$ as in \x{mp}.


By \x{E23} and \x{dtau0}, a direct computation can yield
$$F_{23}(t):= \left.\frac{\partial^2F}{\partial t\partial \delta}\right|_{(\phi_{m,p}(t), t , 0)}=\frac{-8\beta V_0\phi_{m, p}(t)\dot P(t)}{(1-\phi^2_{m,p}(t))^2},$$
and
\begin{align*}
\tau'_{m,p}(0)&=-pT'(\hbar_{m, p})\int^{mT_{v}}_{0}F_{23}(t)\dot \phi_{m,p}(t)dt\\
&=pT'(\hbar_{m, p})\int^{mT_{v}}_{0}\dot P(t)\dot G_{m,p}(t)dt\\
&=pT'(\hbar_{m, p})[\left.\dot P(t)G_{m,p}(t)\right|^{mT_{v}}_{0}-\int^{mT_{v}}_{0}G_{m,p}(t)\ddot P(t)dt],
\end{align*}
where
\begin{equation}\label{E52}
G_{m,p}(t)=\frac{4\beta V_0}{1-\phi^2_{m,p}(t)},
\end{equation}
and
$$\dot G_{m,p}(t)=\frac{8\beta V_0\phi_{m, p}(t)\dot\phi_{m,p}(t)}{(1-\phi^2_{m,p}(t))^2}.$$
Further, from $P(t)=cos\omega_0t$ we have
\begin{equation}\label{E53}
\tau'_{m,p}(0)=\omega^2_0pT'(\hbar_{m, p})\int^{mT_{v}}_{0}G_{m,p}(t)cos\omega_0t~dt.
\end{equation}

From \x{E53}, we see that  $\tau'_\pp(0)=0$ if $ m/(2p)\not \in \N$. This is to say, if $m$ is odd and $1\le p\le \nu_m$, or $m$ is even and $ m/2+1 \le p \le \nu_m$, then the sign of $\tau'_\pp(0)$ depends on $(\pp)$ in a delicate way (see \cite{CC}, Theorem 4.1).

In addition, we consider $ m/(2p)\in\N$, i.e., the following case
   \be \lb{pm09}
    \qq m=2np,
   \ee
where $n, p\in \N$. For simplicity, we denote
    \[
    \phi_{n}(t):= \phi_{2np,p}(t)\equiv\phi_{2n,1}(t).
    \]
By the fact $\hbar_{2np,p}\equiv \hbar_{2n,1}=:\hbar_{n}$ and Theorem \ref{Theo3.1}, it is obtained that $\phi_{n}(t)$ has the minimal period
\be \lb{Th11}
    T_{n}:=T(\hbar_{n})=2nT_{v} \andq   \hbar_{n}= (T_{n})^{-1}(2nT_{v})\in(0, \hbar_\ast].
    \ee
Comparing with $\phi_{m,p}(t)$ for arbitrary integers $m\geq1$ and $p\geq0$, there are more symmetries on $\phi_{n}(t)$ as follows:

\begin{equation}\label{E57}
\left\{
\begin{array}{l}
\phi_{n}(-t)\equiv-\phi_{n}(t),\\
\phi_{n}(t+nT_{v})\equiv-\phi_{n}(t),\\
\phi_{n}(nT_{v}-t)\equiv\phi_{n}(t),\\
\phi_{n}(t)>0,~\mathrm{for}~t\in(0, nT_{v}),\\
\phi_{n}(t)\mbox{ is strictly increasing on $[0,\frac{nT_{v}}{2}]$.}
\end{array}
\right.
\end{equation}
Then, substituting $\phi_{n}(t)$ into \x{E52}, one has
\begin{equation}\label{E58}
G_n(t):=\frac{4\beta V_0}{1-\phi^2_{n}(t)},
\end{equation}
and
\begin{equation}\label{E59}
\left\{
\begin{array}{l}
G_n(t)\mathrm{~is~even~and~has~the~minimal~period}~nT_{v},\\
G_{n}(nT_{v}-t)=G_{n}(t),\\
G_{n}(t)\mbox{ is strictly increasing on $[0,\frac{nT_{v}}{2}]$.}
\end{array}
\right.
\end{equation}
Owing to (\ref{E53}) and (\ref{E59}), we have
\begin{align}\label{E510}
\tau'_{n}(0)&:=\tau'_{2np,p}(0)=\tau'_{2n,1}(0)\notag\\
&=\omega^2_0T'(\hbar_{n})\int^{2nT_{v}}_{0}G_n(t)cos(\omega_0t)dt\notag\\
&=2\omega^2_0T'(\hbar_{n})\int^{nT_{v}}_{0}G_n(t)cos(\omega_0t)dt\notag\\
&=2\omega^2_0T'(\hbar_{n})[\int^{\frac{nT_{v}}{2}}_{0}G_n(t)cos(\omega_0t)dt+\int^{nT_{v}}_{\frac{nT_{v}}{2}}G_n(t)cos(\omega_0t)dt]\notag\\
&=4\omega^2_0T'(\hbar_{n})\int^{\frac{nT_{v}}{2}}_{0}G_n(t)cos(\omega_0t)dt.
\end{align}
Again let $$A_n:=\int^{\frac{nT_{v}}{2}}_{0}G_n(t)cos(\omega_0t)dt.$$
From(\ref{E510}), when $n=1$, we have
\begin{align*}
A_1&=\int^{\frac{T_{v}}{2}}_{0}G_1(t)cos(\omega_0t)dt\notag\\
&=\int^{\frac{T_{v}}{4}}_{0}G_1(t)cos(\omega_0t)dt+\int^{\frac{T_{v}}{2}}_{\frac{T_{v}}{4}}G_1(t)cos(\omega_0t)dt\notag\\
&=\int^{\frac{T_{v}}{4}}_{0}\z(G_1(t)-G_1(\frac{T_{v}}{2}-t)\y)cos(\omega_0t)dt.
\end{align*}
From the third item of (\ref{E59}), we know $G_1(t)$ is strictly increasing on $[0, \frac{T_{v}}{2}]$. Then, we conclude that $A_1<0$, i.e., $\tau'_{1}(0) < 0$. Therefore, $\phi_{1}(t, \delta)$ is elliptic and linearized stable.

Furthermore, we calculate
\begin{align*}
A_{2k-1}=&\int^{\frac{(2k-1)T_{v}}{2}}_{0}G_{2k-1}(t)cos(\omega_0t)dt\notag\\
=&\int^{\frac{T_{v}}{4}}_{0}\{[G_{2k-1}(t)-G_{2k-1}(\frac{T_{v}}{2}-t)]-[(G_{2k-1}(\frac{T_{v}}{2}+t)-G_{2k-1}(T_{v}-t))\notag\\
&-(G_{2k-1}(T_{v}+t)-G_{2k-1}(\frac{3T_{v}}{2}-t))]-\cdots\notag\\
&-[(G_{2k-1}(\frac{(2k-3)T_{v}}{2}+t)-G_{2k-1}((k-1)T_{v}-t))\notag\\
&-(G_{2k-1}((k-1)T_{v}+t)-G_{2k-1}(\frac{(2k-1)T_{v}}{2}-t))]\}cos(\omega_0t)dt,  \qq k=2,3,4, \ldots\notag\\
\end{align*}
and
\begin{align*}
A_{2k}=&\int^{kT_{v}}_{0}G_4(t)cos(\omega_0t)dt\notag\\
=&\int^{\frac{T_{v}}{4}}_{0}\{[(G_{2k}(t)-G_{2k}(\frac{T_{v}}{2}-t))-(G_{2k}(\frac{T_{v}}{2}+t)-G_{2k}(T_{v}-t))]\notag\\
& +[(G_{2k}(T_{v}+t)-G_{2k}(\frac{3T_{v}}{2}-t))-(G_{2k}(\frac{3T_{v}}{2}+t)-G_{2k}(2T_{v}-t))]+\cdots\notag\\
&+[(G_{2k}((k-1)T_{v}+t)-G_{2k}(\frac{(2k-1)T_{v}}{2}-t))\notag\\
&-(G_{2k}(\frac{(2k-1)T_{v}}{2}+t)-G_{2k}(kT_{v}-t))]\}cos(\omega_0t)dt, \qq k=3, 4, \ldots
\end{align*}

Next, we will explore the sign of $A_{2k-1}$ and $A_{k}$. From the expression \x{E58}, we can obtain
the derivative of $G_n(t)$
$$G'_n(t)=8\beta V_0\phi_n\phi'_n(1-\phi^2_n)^{-2}, $$\\
and
\begin{align}\label{G}
G''_n(t)=32\beta V_0(1-\phi^2_n)^{-3}\phi'^2_n\phi^2_n+8\beta V_0(1-\phi^2_n)^{-2}\phi'^2_n+8\beta V_0(1-\phi^2_n)^{-2}\phi_n\phi''_n.
\end{align}
Based on \x{H}, we have $$\phi'^2_n=2H-\phi^2_n+\frac{4\beta V^2_0}{1-\phi^2_n}-4\beta V^2_0.$$
Let
\begin{align}\label{Y1}
Y_n=1-\phi^2_n.
\end{align}
Then, one has  $\phi_n=(1-Y_n)^{\frac{1}{2}}$ and
\begin{align}\label{S11}
\phi'^2_n=2H-(1-Y_n)+\frac{4\beta V^2_0}{Y_n}-4\beta V^2_0.
\end{align}
Again from \x{H}, we have $$Y_n(0)=1-\phi^2_n(0)=1$$ and
\begin{align*}
y_1&:=Y_n(\frac{nT_{v}}{2})=1-\phi^2_n(\frac{nT_{v}}{2})\notag\\
&=\frac{1}{2}-H+2\beta V^2_0-\sqrt{4\beta^2V^4_0-4\beta V^2_0H-2\beta V^2_0+H^2+\frac{1}{4}-H}
\end{align*}
Obviously, when $t\in[0, \frac{nT_{v}}{2}]$, we see that $Y_n(t)\in[y_1, 1]$ and $Y_n(t)$ decreases with respect to $t$.

Moreover, from (\ref{E31}), we have
\begin{align}\label{S12}
\phi''_n=-(1-Y_n)^{\frac{1}{2}}(1-\frac{4\beta V^2_0}{Y^2_n})
\end{align}
Substitute \x{Y1}, \x{S11} and \x{S12} into \x{G}, we derive
\begin{align}\label{S13}
G''_n(t)=8\beta V_0Y^{-4}_n[-2Y^{3}_n+(12\beta V^2_0+6-6H)Y^2_n+(8H-4-32\beta V^2_0)Y_n+20\beta V^2_0],
\end{align}
where $Y_n\in[y_1, 1]$. Let $$U_n=-2Y^{3}_n+(12\beta V^2_0+6-6H)Y^2_n+(8H-4-32\beta V^2_0)Y_n+20\beta V^2_0.$$

By MATHEMATICA software, we obtain that there are only one real zero point $y_{0}$ of $U_{n}$ when $Y_{n}\in \mathbb{R}$.  It is note that $U_{n}$ is continuous on $Y_n\in\mathbb{R}$ and $U_n\mid_{Y_n=0}=20\beta V^2_0>0$, $U_n\mid_{Y_n=1}=2H>0$. Hence, whether the real zero point $y_{0}\in[0,1]$ or $y_{0}\notin[0,1]$,
we can conclude $U_{n}\geq 0$ when $Y_{n}\in[0,1]$. Further, one has $U_{n}\geq 0$ when $Y_{n}\in[y_{1},1]$, which imply that $G''_n(t)\geq0$ when $t\in[0, \frac{nT_{v}}{2}]$ and $G_n(t)$ is convex function on $[0, \frac{nT_{v}}{2}]$. Then, from the convexity and monotonicity of $G_n(t)$, we obtain that $A_{2k-1}<0$ and $A_{2k}>0$, respectively. So, we obtain the following theorem.

\bb{Theorem} \lb{M07}
Assuming that $0<\omega_{0}<2n\sqrt{1-4\beta V^{2}_{0}}$,
for $0 < \delta\ll \Delta_0$, any $ p, \ n\in \N$ and $m=2np$, we have

   {\rm (i)}  when $n$ is odd, $\tau'_{n}(0)<0$, i.e., odd $(m,p)$-periodic solutions $\phi_{m,p}(t,\delta)$ of the Comb-Drive finger actuator model \x{E21} are linearly stable.

   {\rm (ii)}  when $n$ is even, $\tau'_{n}(0)>0$, i.e., odd $(m,p)$-periodic solutions $\phi_{m,p}(t,\delta)$ of the Comb-Drive finger actuator model \x{E21} are hyperbolic and Lyapunov unstable.
\end{Theorem}


\subsection{Analytical results for stability of even periodic orbits}

In this subsection, we will analyze the stability of of the family $\varphi_{m,p}(t, \delta)$ of even $(m,p)-$periodic solutions of the Eq. \x{E23} for $m$ and $p$ be as in (\ref{mp}).

By \x{E23} and \x{dtau0}, we have
$$\hat{F}_{23}(t)=\left.\frac{\partial^2F}{\partial t\partial \delta}\right|_{(\varphi_{m,p}(t), t , 0)}=\frac{-8\beta V_0\varphi_{m, p}(t)P'(t)}{(1-\varphi^2_{m,p}(t))^2},$$
\begin{align}\label{E512}
\hat{\tau}'_{m,p}(0)&=-pT'(\hbar_{m, p})\int^{mT_{v}}_{0}\hat{F}_{23}(t)\dot \varphi_{m,p}(t)dt\notag\\
&=pT'(\hbar_{m, p})\int^{mT_{v}}_{0}\dot P(t)\hat{G}'_{m,p}(t)dt\notag\\
&=\omega^2_0pT'(\hbar_{m, p})\int^{mT_{v}}_{0}\hat{G}_{m,p}(t)cos\omega_0t~dt
\end{align}
where
\begin{equation}\label{E513}
\hat{G}_{m,p}(t)=\frac{4\beta V_0}{1-\varphi^2_{m,p}(t)}.
\end{equation}

When $m=2np$, we note that $\phi_{2np,p}(t)$ and $\varphi_{2np,p}(t)$ have the same energy $\hbar_{2np,p}=\hbar_{n}$ and the same minimal period $2nT_{v}$. Let $\varphi_{n}(t):=\varphi_{2np,p}(t)=\phi_{n}(t+\frac{nT_{v}}{2})$. Using the notations in (\ref{E52}) and (\ref{E513}), we have
$$\hat{G}_{n}(t):=\hat{G}_{2np,p}(t)=\hat{G}_{2n,1}(t)=G_{n}(t+\frac{nT_{v}}{2}).$$
Hence,
\begin{align*}
\hat{\tau}'_{n}(0)&:=\hat{\tau}'_{2np,n}(0)=\hat{\tau}'_{2n,1}(0)\\
&=\omega^2_0T'(\hbar_{n})\int^{2nT_{v}}_{0}\hat{G}_{n}(t)cos\omega_0t~dt\\
&=\omega^2_0T'(\hbar_{n})\int^{2nT_{v}}_{0}G_{n}(t+\frac{nT_{v}}{2})cos\omega_0t~dt\\
&=\omega^2_0T'(\hbar_{n})\int^{2nT_{v}+\frac{nT_{v}}{2}}_{\frac{nT_{v}}{2}}G_{n}(t)cos\omega_0(t-\frac{nT_{v}}{2})~dt\\
&=(-1)^n\omega^2_0T'(\hbar_{n})\int^{2nT_{v}}_{0}G_{n}(t)cos\omega_0t~dt\\
&=(-1)^n\tau'_{n}(0).
\end{align*}

Based on the above relation, we have the following theorem.
\bb{Theorem} \lb{M5}
Assuming that $0<\omega_{0}<2n\sqrt{1-4\beta V^{2}_{0}}$, for any $ p, \ n\in \N$, $m=2np$, we have
    \[
   \hat\tau'_{n}(0)>0.
   \]
Consequently, for $0 < \delta\ll \Delta_0$, even $(m,p)$-periodic solutions $\varphi_{m,p}(t,e)$ of the comb-drive finger actuator model \x{E21} are hyperbolic and Lyapunov unstable.
   \end{Theorem}

\section{Conclusions}\label{sec:4}

In this paper, we analytically studied the linear stability/instability of the nonconstant symmetric periodic solutions of the comb-drive finger actuator model. In the light of the relationship between the potential and the period as a function of the energy, we firstly deduced the three properties of the minimal period $T(\hbar)$ of the solution of the corresponding autonomous equation. For the prescribed energy range $ \hbar\in(0, \hbar_\ast)$, the minimal period $T(\hbar)$ satisfies $T(\hbar)\in(2\pi/\sqrt{1-4\beta V^2_0},+\oo)$ and $T'(\hbar)>0$, implying that $T(\hbar)$ is a strictly monotone increasing function on energy $\hbar$. Then, based on $T'(\hbar)>0$ and the stability criteria, we arrived at the linear stability/instability of the symmetric $(m,p)$-periodic solutions when the parameter $\delta$ is small. Specifically, for $m=2np$, odd $(m,p)$-periodic solutions $\phi_{m,p}(t,e)$ are linearly stable when $n$ is odd, and $\phi_{m,p}(t,e)$ are hyperbolic and Lyapunov unstable when $n$ is even.  However, all of the even $(m,p)$-periodic solutions  $\vp_{m,p}(t,e)$ are hyperbolic and Lyapunov unstable.

\end{document}